\documentclass[12pt]{amsart}
\usepackage{geometry} % see geometry.pdf on how to lay out the page. There's lots.
\newtheorem{theorem}{Theorem}

\title{Efficient computation of tight approximations to Chernoff bounds}
\author{D. K. L. Shiu}
\begin{document}

\maketitle
%\tableofcontents

\section{Introduction}
Chernoff bounds are a powerful application of the Markov inequality to produce strong bounds on the tails of probability distributions. They are often used to bound the tail probabilities of sums of Poisson trials, or in regression to produce conservative confidence intervals for the parameters of such trials. The bounds provide expressions for the tail probabilities that can be inverted for a given probability/confidence to provide tail intervals. The inversions involve the solution of transcendental equations and it is often convenient to substitute approximations that can be exactly solved e.g. by the quadratic equation.

In this paper we introduce approximations for the Chernoff bounds whose inversion can be exactly solved with a quadratic equation, but which are closer approximations than those adopted previously.

\section{Approximating tail distributions}
We begin with the predictive application for sums of independent Poisson trials $X_i$ where the sum $\sum X_i$ has expectation $\mu$. We have (see for example \cite{mitzenmacher}, Theorem 4.4 and 4.5) the Chernoff bound tail probabilities:
$$\mathbb P\left(\sum_iX_i\ge(1+\delta)\mu\right)\le\exp\left(\left(\delta-(1+\delta)\log(1+\delta)\right)\mu\right)$$
$$\mathbb P\left(\sum_iX_i\le(1-\delta)\mu\right)\le\exp\left(\left(-\delta-(1-\delta)\log(1-\delta)\right)\mu\right).$$

If we are given an upper bound tail probability, $\gamma$ and we might be required to find $\delta_U$ such that $\mathbb P(\sum X_i\ge (1+\delta_U)\mu)<\gamma$ or such that
$\mathbb P(\sum X_i\le (1-\delta_L)\mu)<\gamma$. We could solve
$$\delta_U-(1+\delta_U)\log(1+\delta_U)=\frac{\log\gamma}\mu$$
or
$$-\delta_L-(1-\delta_L)\log(1-\delta_L)=\frac{\log\gamma}\mu$$
numerically (e.g. by iterative methods such as binary search or Newton's method). 

In practice approximations $\log(1-\delta)\ge-\delta-\delta^2/2$ and $\log(1+\delta)\ge2\delta/(2+\delta)$ valid for $0\le\delta<1$ are used to provide the more wieldy bounds
\begin{equation}\mathbb P\left(\sum_iX_i\ge(1+\delta)\mu\right)\le\exp\left(\frac{-\delta^2\mu}{2+\delta}\right)\label{eq:old_u}\end{equation}
\begin{equation}\mathbb P\left(\sum_iX_i\le(1-\delta)\mu\right)\le\exp\left(\frac{-\delta^2\mu}2\right).\label{eq:old_l}\end{equation}
For our tail probability $\gamma$ we can use these looser expressions to form the quadratic equations
$$\delta_U^2+\left(\frac{\log\gamma}\mu\right)\delta_U+2\left(\frac{\log\gamma}\mu\right)=0$$
$$\delta_L^2+2\left(\frac{\log\gamma}\mu\right)=0$$
and evaluate the closed form expressions
$$\delta_U=\frac{-\frac{\log\gamma}\mu+\sqrt{\left(\frac{\log\gamma}\mu\right)^2-8\frac{\log\gamma}\mu}}2$$
$$\delta_L=\sqrt{-\frac{2\log\gamma}\mu}$$
and conclude that $\mathbb P(\sum X_i\ge (1+\delta_U)\mu)<\gamma$ and $\mathbb P(\sum X_i\le (1-\delta_L)\mu)<\gamma$. If a two-tailed bound is required, we can note that $\delta_U\ge\delta_L$ so that $\mathbb P((1-\delta_U)\mu\le\sum X_i\le (1+\delta_U)\mu)<1-2\gamma$.

The inequalities $\log(1-\delta)\ge-\delta-\delta^2/2$ and $\log(1+\delta)\ge2\delta/(2+\delta)$ are quadratic Pad\'e approximations to $\log(1\pm\delta)$, but we can do better with approximations to the functions
\begin{equation}\delta-(1+\delta)\log(1+\delta)=-\frac{\delta^2}2+\frac{\delta^3}6-\frac{\delta^4}{12}+\frac{\delta^5}{20}-\frac{\delta^6}{30}-\cdots\label{eq:power_u}\end{equation}
and
\begin{equation}-\delta-(1-\delta)\log(1-\delta)=-\frac{\delta^2}2-\frac{\delta^3}6-\frac{\delta^4}{12}-\frac{\delta^5}{20}-\frac{\delta^6}{30}-\cdots\label{eq:power_l}\end{equation}
We see that
\begin{equation}\frac{-3\delta^2}{6+2\delta}=-\frac{\delta^2}2+\frac{\delta^3}6-\frac{\delta^4}{18}+\frac{\delta^5}{54}-\frac{\delta^6}{162}-\cdots\label{eq:pade_u}\end{equation}
and
\begin{equation}\frac{-9\delta^2}{18-6\delta-\delta^2}=-\frac{\delta^2}2-\frac{\delta^3}6-\frac{\delta^4}{12}-\frac{\delta^5}{27}-\frac{11\delta^6}{648}-\cdots\label{eq:pade_l}\end{equation}
It is easy to confirm that for $\delta\in(0,1)$ we have
$$\delta-(1+\delta)\log(1+\delta)<\frac{-3\delta^2}{6+2\delta}$$
and
$$-\delta-(1-\delta)\log(1-\delta)<\frac{-9\delta^2}{18-6\delta-\delta^2},$$
so that
\begin{equation}\mathbb P\left(\sum_iX_i\ge(1+\delta)\mu\right)\le\exp\left(\frac{-3\delta^2}{6+2\delta}\right)\label{eq:bernstein}\end{equation}
$$\mathbb P\left(\sum_iX_i\le(1-\delta)\mu\right)\le\exp\left(\frac{-9\delta^2}{18-6\delta-\delta^2}\right).$$
Note that the difference between the power series in (\ref{eq:power_u}) and (\ref{eq:pade_u}) is $O(\delta^4)$ and the difference between the power series in (\ref{eq:power_l}) and (\ref{eq:pade_l}) is $O(\delta^5)$. For comparison, the bounds in (\ref{eq:old_u}) and (\ref{eq:old_l}) introduce a $O(\delta^3)$ difference into the exponent. We also see that if we restrict to identically distributed trials (i.e.Bernoulli rather than Poisson), equation (\ref{eq:bernstein}) can be rearranged to an estimate previously derived using Bernstein inequalities \cite{bernstein}. For our tail probability $\gamma$ we can use these more accurate expressions to form the quadratic equations
$$3\delta_U^2+2\left(\frac{\log\gamma}\mu\right)\delta_U+6\left(\frac{\log\gamma}\mu\right)=0$$
$$\left(9-\frac{\log\gamma}\mu\right)\delta_L^2-6\left(\frac{\log\gamma}\mu\right)\delta_L+18\left(\frac{\log\gamma}\mu\right)=0$$
and evaluate the closed form expressions
$$\delta_U=\frac{-\left(\frac{\log\gamma}\mu\right)+\sqrt{\left(\frac{\log\gamma}\mu\right)^2-18\left(\frac{\log\gamma}\mu\right)}}3$$
$$\delta_L=\frac{3\left(\frac{\log\gamma}\mu\right)+\sqrt{9\left(\frac{\log\gamma}\mu\right)^2-18\left(\frac{\log\gamma}\mu\right)\left(9-\frac{\log\gamma}\mu\right)}}{9-\frac{\log\gamma}\mu}.$$

We have now proven the following:
\begin{theorem} Let $X_1,\ldots, X_n$ be independent Poisson trials. Let $X=\sum_{i=1}^nX_i$ and $\mu=\mathbb E[X]$. Let $0<\gamma<1$ be a fixed tail probability and write $\beta=(\log\gamma)/\mu$. Let
$$\delta_U=\frac{-\beta+\sqrt{\beta^2-18\beta}}3$$
$$\delta_L=3\left(\frac{\beta+\sqrt{\beta^2-2\beta\left(9-\beta\right)}}{9-\beta}\right).$$
Then
$$\mathbb P(X\ge (1+\delta_U)\mu)<\gamma$$
and
$$\mathbb P(X\le (1-\delta_L)\mu)<\gamma.$$
\end{theorem}

\section{Approximating confidence intervals}
For our regression application, we have an observed value $\hat\mu$ of $\sum X_i$ and a target confidence $\gamma$.  Our goal is to identify a range of possible underlying $\mu$ values such that the likelihood of our observation $\hat\mu$ for $\mu$ in that range is less than or equal to $\gamma$. The complement of the range then provides a conservative confidence interval for $\mu$ with confidence at least $1-\gamma$. 

We first develop Chernoff bounds of a slightly different form. Starting from the generic Markov bound applied to $e^{tX}$ we have
$$\mathbb P(X\le a)\le\min_{t<0} \frac{\mathbb E[e^{tX}]}{e^{ta}},\qquad \mathbb P(X\le a)\ge\min_{t>0} \frac{\mathbb E[e^{tX}]}{e^{ta}}.$$
Following the usual argument, we note that $\mathbb E[e^{tX_i}]=1+p_i(e^t-1)\le\exp(p_i(e^t-1))$ and so $\mathbb E[e^{tX}]\le\exp((e^t-1)\mu)$. Then with $t=-\log(1+\delta)$ and $t=-\log(1-\delta)$ we have
$$\mathbb P\left(\sum_iX_i\le\frac\mu{(1+\delta)}\right)\le\exp\left(\left(-\frac{\delta}{(1+\delta)}+\frac{\log(1+\delta)}{(1+\delta)}\right)\mu\right)$$
and
$$\mathbb P\left(\sum_iX_i\ge\frac\mu{(1-\delta)}\right)\le\exp\left(\left(\frac{\delta}{(1-\delta)}+\frac{\log(1-\delta)}{(1-\delta)}\right)\mu\right).$$

Using the above bounds, if we let $\mu_\ell$ be the value for which $\hat\mu=\mu_\ell/(1+\delta_\ell)$ where $\delta_\ell$ is a solution to
\begin{equation}\exp\left(\left(-\delta_\ell+\log(1+\delta_\ell)\right)\hat\mu\right)=\gamma\label{eq:exactlower}\end{equation}
then by monotonicity, the interval $(-\infty,\mu_\ell]$ is a suitable range of exceptional $\mu$ values and $(\mu_\ell,\infty)$ is a suitable conservative confidence interval. For a conservative confidence interval that gives an upper bound for $\mu$, we can, by a similar process, find the $\mu_u=\hat\mu(1-\delta_u)$ where
\begin{equation}\exp\left(\left(\delta_u+\log(1-\delta_u)\right)\hat\mu\right)=\gamma\label{eq:exactupper}\end{equation}
and develop the conservative confidence interval $(-\infty,\mu_u)$. We can even combine our calculation for the conservative confidence interval $(\mu_\ell,\mu_u)$ in which we would have confidence at least $1-2\gamma$.

As in the previous section, the equations (\ref{eq:exactupper}) and (\ref{eq:exactlower}) can be solved numerically by iterative methods such as binary search or Newton's method. However, we seek an expression that can be solved using the quadratic formula.

For a lower bound for $\mu$ with level of confidence $1-\gamma$ with $0<\gamma<1$, we therefore aim to identify the values $\delta$ such that
\begin{equation}(\delta+\log(1-\delta))\hat\mu\le\log\gamma.\label{eq:multlower}\end{equation}
By monotonicity, for $0<\delta<1$, identifying the value where equality is attained proves the bound for all greater values. For a similarly confident upper bound we need to identify the values $\delta$ such that
\begin{equation}(-\delta+\log(1+\delta))\hat\mu\le\log\gamma\label{eq:multupper}\end{equation}
and again, identifying the value where equality is attained proves the bound for all greater values

For $0<\delta<1$ we have the Pad\'e approximation to $-\delta+\log(1+\delta)$
\[-\delta+\log(1+\delta)<\frac{-3\delta^2}{6+4\delta}.\label{eq:padeupper}\] 
and so if we take $\delta_U$ to be the positive root of
\[3\delta_U^2+4\frac{\log\gamma}{\hat\mu}\delta_U+6\frac{\log\gamma}{\hat\mu}\]
then we have
$$(-\delta_U+\log(1+\delta_U))\hat\mu\le\log\gamma.$$
It follows that $\mu\ge(1+\delta_U)\hat\mu$ with probability at most $\gamma$ so that $\mu<(1+\delta_U)\hat\mu$ with probability at least $1-\gamma$. We note the power series expansions for $0<\delta<1$
$$-\delta+\log(1+\delta)=-\frac{\delta^2}2+\frac{\delta^3}3-\frac{\delta^4}4+\frac{\delta^5}5-\frac{\delta^6}6-\cdots$$
$$\frac{-3\delta^2}{6+4\delta}=-\frac{\delta^2}2+\frac{\delta^3}3-\frac{2\delta^4}9+\frac{4\delta^5}{27}-\frac{8\delta^6}{81}-\cdots$$
so that our approximation differs by $O(\delta^4)$.

Likewise we also have the quadratic Pad\'e approximation to $\delta+\log(1-\delta)$
\begin{equation}\delta+\log(1-\delta)<\frac{-9\delta^2}{18-12\delta-\delta^2}\label{eq:padelower}\end{equation}
From (\ref{eq:padelower}) we conclude that if $\delta_L$ is taken to be the positive root of the equation
\[\left(9-\frac{\log\gamma}{\hat\mu}\right)\delta_L^2-12\frac{\log\gamma}{\hat\mu}\delta_L+18\frac{\log\gamma}{\hat\mu}=0\]
then we have
$$(\delta_L+\log(1-\delta_L))\hat\mu<\log\gamma.$$
It follows that $\mu\le(1-\delta_L)\hat\mu$ with probability at most $\gamma$ so that $\mu>(1-\delta_L)\hat\mu$ with probability at least $1-\gamma$.  Again by considering the power series expansions for $0<\delta<1$
$$\delta+\log(1-\delta)=-\frac{\delta^2}2-\frac{\delta^3}3-\frac{\delta^4}4-\frac{\delta^5}5-\frac{\delta^6}6-\cdots$$
$$\frac{-9\delta^2}{18-12\delta-\delta^2}=-\frac{\delta^2}2-\frac{\delta^3}3-\frac{\delta^4}4-\frac{5\delta^5}{27}-\frac{89\delta^6}{648}-\cdots$$
so that our approximation differs by $O(\delta^5)$.

The two estimates can be combined into a two-ended confidence interval allowing us to conclude that $\mu\in((1-\delta_L)\hat\mu,(1+\delta_U)\hat\mu)$ with probability at least $1-2\gamma$. Alternatively, if a symmetric expression is desired, we note that for $0<\delta<1$ we have
\[-\delta+\log(+\delta)<-\frac{3\delta^2}{6+4\delta} <\delta+\log(1-\delta)<-\frac{9\delta^2}{18-12\delta-\delta^2}\]
and so by monotonicity $\delta_U>\delta_L$ and the interval $\mu\in((1-\delta_U)\hat\mu,(1+\delta_U)\hat\mu)$ can be used with confidence at least $1-2\gamma$. 

\begin{theorem} Let $X_1,\ldots, X_n$ be independent Poisson trials. Let $X=\sum_{i=1}^nX_i$ and suppose that we have a sample $\hat\mu$ from $X$ . Let $0<\gamma<1$ be a fixed bound on confidence and write $\beta=(\log\gamma)/\hat\mu$. Let
$$\delta_U=\frac{-2\beta+\sqrt{4\beta^2-18\beta}}{3}$$
$$\delta_L=\frac{6\beta+\sqrt{36\beta^2-18\beta(9-\beta)}}{9-\beta}.$$
Then with confidence at least $1-\gamma$ we can say
$$\mathbb E[X]<(1+\delta_U)\hat\mu$$
and similarly with confidence at least $1-\gamma$ we can say 
$$\mathbb E[X]>(1-\delta_L)\hat\mu.$$
\end{theorem}

\section{Numerical examples}
We consider examples using Bernoulli trials which are a frequent use of such bounds. Suppose that we have $\mathbb P(X_i=1)=0.0002$ and that we run 1,000,000 trials. We have $\mu=200$. We consider tail probabilities $\gamma$ of 0.05, 0.01, 0.000000002 (corresponding to ``six sigma''), and 5.421e-20 (corresponding to a probability of $2^{-64}$ which is relevant to failure rates in cryptography). For each $\gamma$ we compute $\delta_U$ and $\delta_L$ using the exact transcendental Chernoff formulae, the old quadratic formulae and the new quadratic formulae of this paper. Solutions are given to four significant figures

\medskip

\begin{tabular}{|c|cc|cc|cc|}
\hline
$\gamma$ & Exact $\delta_U$ & Exact $\delta_L$ & Old $\delta_U$ & Old $\delta_L$ & New $\delta_U$ & New $\delta_L$ \\
\hline
0.05               & 0.1780 & 0.1680 & 0.1807 & 0.1731 & 0.1781 & 0.1680 \\ 
0.01               & 0.2221 & 0.2068 & 0.2264 & 0.2146 & 0.2224 & 0.2068 \\
0.000000002 & 0.4798 & 0.4127 & 0.5004 & 0.4476 & 0.4822 & 0.4133 \\
5.421e-20      & 0.7365 & 0.5870 & 0.7861 & 0.6660 & 0.7441 & 0.5898 \\
\hline
\end{tabular}

\medskip

As we expected, our new approximation is closer to the exact Chernoff bound, particularly for larger $\delta$ which correspond to smaller $\gamma$ or smaller $\mu$. For the smallest $\gamma$ value, we observe that the upper bound for the number of successful trials is 347, 357, and 348 respectively and that the lower bound is 83, 67, and 83 respectively.

Turning now to our regression estimates, we assume that we run 1,000,000 Bernoulli trials with unknown probability and that we observe $\hat\mu=212$ successes. Again we choose confidence levels $1-\gamma$ with $\gamma=0.05, 0.01, 0.000000002, 5.421e-20$. For each $\gamma$ we compute $\delta_U$ and $\delta_L$ using the exact transcendental Chernoff formulae and the new quadratic formulae of this paper. Solutions are given to four significant figures

\medskip

\begin{tabular}{|c|cc|cc|}
\hline
$\gamma$ & Exact $\delta_U$ & Exact $\delta_L$ & New $\delta_U$ & New $\delta_L$ \\
\hline
0.05               & 0.1777 & 0.1588 & 0.1778 & 0.1588  \\ 
0.01               & 0.2232 & 0.1942 & 0.2234 & 0.1942  \\
0.000000002 & 0.4998 & 0.3741 & 0.5022 & 0.3746  \\
5.421e-20      & 0.7933 & 0.5156 & 0.8013 & 0.5176  \\
\hline
\end{tabular}

\section{Higher degree approximation}
Eager readers will be aware that cubic and quartic equation also admit closed form solutions. For still greater accuracy, higher degree Pad\'e approximations could be used. We note the following Pad\'e approximations that could be used to this end.

$$x-(1+x)\log(1+x)<\frac{-15x^2-7x^3}{30+24x+3x^2}$$

$$x-(1+x)\log(1+x)<\frac{-210x^2-200x^3-35x^4}{420+540x+180x^2+12x^3}$$

$$-x-(1-x)\log(1-x)<\frac{-210x^2+125x^3}{420-390x+60x^2+3x^3}$$

$$-x-(1-x)\log(1-x)<\frac{7350x^2-8260x^3+1975x^4}{-14700+21420x-8640x^2+780x^3+18x^4}$$

$$-x+\log(1+x)<\frac{240x^2+155x^3}{-480-630x-180x^2+3x^3}$$

$$-x+\log(1+x)<\frac{-210x^2-220x^3-45x^4}{420+720x+360x^2+48x^3}$$

$$x+\log(1-x)<\frac{-240x^2+155x^3}{480=630x+180x^2+x^3}$$

$$x+\log(1-x)<\frac{3150x^2-3780x^3+985x^4}{-6300+11760-6660x^2+1080x^3+6x^4}$$

\end{document}